\documentclass[11pt]{article}
\usepackage{graphicx}
\usepackage{color}
\usepackage{amssymb}

\begin{document}
\newtheorem{defn0}{Definition}[section]
\newtheorem{prop0}[defn0]{Proposition}
\newtheorem{thm0}[defn0]{Theorem}
\newtheorem{lemma0}[defn0]{Lemma}
\newtheorem{coro0}[defn0]{Corollary}
\newtheorem{exa}[defn0]{Example}
\def\rig#1{\smash{ \mathop{\longrightarrow}
    \limits^{#1}}}
\def\nwar#1{\nwarrow
   \rlap{$\vcenter{\hbox{$\scriptstyle#1$}}$}}
\def\near#1{\nearrow
   \rlap{$\vcenter{\hbox{$\scriptstyle#1$}}$}}
\def\sear#1{\searrow
   \rlap{$\vcenter{\hbox{$\scriptstyle#1$}}$}}
\def\swar#1{\swarrow
   \rlap{$\vcenter{\hbox{$\scriptstyle#1$}}$}}
\def\dow#1{\Big\downarrow
   \rlap{$\vcenter{\hbox{$\scriptstyle#1$}}$}}
\def\up#1{\Big\uparrow
   \rlap{$\vcenter{\hbox{$\scriptstyle#1$}}$}}
\def\lef#1{\smash{ \mathop{\longleftarrow}
    \limits^{#1}}}
\def\O{{\cal O}}
\def\L{{\cal L}}
\def\M{{\cal M}}
\def\N{{\cal N}}
\newcommand{\R}{\mathbb{R}}
\newcommand{\C}{\mathbb{C}}
\newcommand{\K}{\mathbb{K}}
\def\P#1{{\bf P}^#1}
\newcommand{\defref}[1]{Def.~\ref{#1}}
\newcommand{\propref}[1]{Prop.~\ref{#1}}
\newcommand{\thmref}[1]{Thm.~\ref{#1}}
\newcommand{\lemref}[1]{Lemma~\ref{#1}}
\newcommand{\corref}[1]{Cor.~\ref{#1}}
\newcommand{\exref}[1]{Example~\ref{#1}}
\newcommand{\secref}[1]{Section~\ref{#1}}

\newcommand{\qedd}{\hfill\framebox[2mm]{\ }\medskip}
\newcommand{\codim}{\textrm{codim}}
\newtheorem{rem0}[defn0]{Remark}
\def\proof{{\it Proof:\  \ }}
\definecolor{purple}{rgb}{0.6,0,0.3}

\author{Pierre Comon, Giorgio Ottaviani\footnote{The second author is member of GNSAGA-INDAM.}}
\title{On the typical rank of real binary forms}
\date{}
\maketitle
%\tableofcontents

\begin{abstract}
We determine the rank of a general real binary form of degree $d=4$ and $d=5$.
In the case $d=5$, the possible values of the rank of such general forms are $3, 4, 5$.
The existence of three typical ranks was unexpected. We prove that a real binary form of degree $d$
with $d$ real roots has rank $d$.
\end{abstract}

\bigskip

\section{Introduction}

\textbf{Motivation.}
The decomposition of tensors into a minimal sum of rank-1 terms, sometimes referred to as the \emph{Canonical decomposition} (CanD), has raised a great interest in several quite different branches of engineering \cite{SidiBG00:ieeesp,Kroo08,SmilBG04,BadeK06:tms,MoruHHPA06:ni,CarrC70:psy,Como02:oxford} as well as in mathematics \cite{HackKT05:jnm,IarrK99,BurgCS97,Lick85:laa,lanteit,BGI,BL}. More precisely, if we take the example of tensors of order 4 to fix the ideas, the CanD of a tensor $T$ is defined as:
$$
T_{ijk\ell} = \sum_{p=1}^r \alpha_i(p)\,\beta_j(p)\,\gamma_k(p)\,\delta_\ell(p)
$$
The minimal value of $r$ such that the equality holds is called the \emph{rank} of tensor $T$.

In applications, these tensors are often corrupted by measurement or estimation noises, which can be modeled as additive random variables with continuous probability distribution. The consequence is that tensors to be decomposed are \emph{generic}. A key question is then to know what is the rank of such tensors as a function of order and dimensions, and when they admit a \emph{finite number} of CanD. 
For real or complex \emph{symmetric} tensors decomposed in the \emph{complex} field, this question has been answered, and enters the framework of the \emph{Waring problem} \cite{IarrK99,cilisurvey}. The Theorem of Alexander and Hirschowitz (see \cite{AH} and also \cite{BO})
answers that the generic rank is always the expected one, with a finite list of exceptions.
However, it has not yet received any response either for non symmetric tensors, or for decompositions in the real field.
ten Berge showed that real third order generic tensors can assume two typical ranks, $r$ and $r+1$ \cite{Tenb00:psy}. But he could not find generic tensors assuming more than two typical ranks.
Regarding general complex tensors, some partial results are contained in \cite{AOP}, while \cite{CGG} provides a complete answer for
tensors of format $2\times\ldots\times 2$.
We point out also the recent preprint \cite{Rez}, where the related concept of signature is studied.

\textbf{Contribution.}
In order to exhibit more than two typical ranks, it seems necessary to consider tensors of order higher than 3. The simplest case we can think of is that of real symmetric tensors of dimension 2 and order 4 or 5. Yet, such tensors are bijectively associated with homogeneous polynomials in two variables \cite{EhreR93:ejc}, namely binary quartics and quintics, respectively.

The previous definition translates to binary forms as follows. The rank of a binary form $f(x,y)=\sum_{i=0}^d{d\choose i}\,a_i\,x^{d-i}y^i$ of degree $d$
defined over the field $\K$
is the minimum integer $r$ such that there is a decomposition $f=\sum_{j=1}^rc_j\,(l_j)^d$
where $l_j$ are linear forms and $c_j\in \K$. 
Note that the rank is called ``length'' in \cite{Rezn92:mams}.
If $\K$ is the real field, coefficients $c_j$ can be imposed to belong to $\{-1,\,1\}$. 
If $\K$ is the complex field, one can impose $c_j=1$ for all coefficients. 
In the complex field, the  rank of a general binary form $f$ of odd degree $d=2n+1$ is $n+1$. The Sylvester Theorem asserts that
 the decomposition of such general form $f$ as a sum of $n+1$ powers of linear forms is unique (up to a multiplicative factor for each $c_j$) , and gives also a way to determine it. 
 The  {\it rank} of a general binary form $f$ of even degree $d=2n$ is $n+1$, but in this case such decompositions form an infinite set, 
which can be identified with the projective line.
 
 It is interesting that no semicontinuity holds for the rank.
Indeed there are forms such that their rank is bigger than the general value
 and other such that the rank is lower than the general value.
 Indeed it is well known (see for example Comas and Seiguer \cite{ComasS01:arxiv}) that on the complex numbers the rank
 can attain all the values between $1$ and $d$.
 
 Define $S_{d,r}^\K=\{f\in Sym^d(\K^2)|\textrm{rank\ }f=r\}$, the set of binary forms of rank exactly $r$. On complex numbers, $S_{d,r}^{\C}$ has a non empty interior only for the pairs $(d,r)$ where $r=\lfloor\frac{d}{2}\rfloor+1$. That is, for every $d$, there is a unique $r$ such that $S_{d,r}^{\C}$ is dense. The latter values of $r$ are called \emph{generic ranks} in $\C$.
 On the real numbers we are intested in the pairs $(d,r)$ such
that $S_{d,r}^{\R}$ has non empty interior. If for a given value of $d$ there exist several values of $r$, the latter are called \emph{typical ranks}. 
In other words, if $r_i$ are typical, none of the sets $S_{d,r_i}^{\R}$ is dense, but $\cup_{i}S_{d,r_i}^{\R}$ is.
For $d=2, 3$ the situation
is well known and we recall it in the next section.
Our main results are the following, for $d=4, 5$.

\medskip

\textbf{Main Theorem}\\
 (i) $S_{4,r}^{\R}$ has non empty interior
 only for $r=3,4$.
 More precisely, assume that $f$ has distinct roots, then
 \begin{itemize}
 \item{(i1)} if $f$ has four  real roots then the rank of $f$ is $4$.
 \item{(i2)} if $f$ has zero or two real roots then the rank of $f$ is $3$.
   \end{itemize}
 (ii) $S_{5,r}$ has non empty interior
 only for $r=3,4, 5$. More precisely, assume that $f$ has distinct roots, then
 \begin{itemize}
 \item{(i1)} if $f$ has five  real roots then the rank of $f$ is $5$.
 \item{(i2)} if $f$ has one or three real roots then the rank of $f$ is $3$ or $4$, according to the sign of an invariant of degree $12$
defined in section $4$, where we provide an algorithm and some pictures).
\end{itemize}

 \bigskip
 
Very probably, $S_{d,r}$ has a non empty interior if and only if $\lfloor\frac{d}{2}\rfloor+1\le r\le d$.
 We leave this as a conjecture to the readers (proved here for $d\le 5$), as the intriguing problem of the relations between $S_{d,r}$ and the loci
where the number of real roots is fixed.
We prove also (see the \propref{dreal}) that if a form of degree $d$ has the maximal number $d$ of real roots,
then its rank attains the maximum possible value which is again $d$.

 There is always an open subset where the general rank is the same as the complex one. In other words, for given order and dimension, the smallest typical rank in the real field coincides with the generic rank in the complex field (see \cite{ComoTDC09:laa},
\cite{Fr}).

\section{Preliminaries: quadratic and cubic forms}

An elementary example illustrating the difference between real and complex rank is given by the equality
$$2x^3-6xy^2=(x+\sqrt{-1}y)^3+(x-\sqrt{-1}y)^3=(2x)^3-(x+y)^3-(x-y)^3$$
In this case the complex rank is $2$ and the real rank is $3$.

Let's begin with a known general statement.
\begin{prop0}\label{maxrank}
Any binary real form of degree $d$ has real rank $\le d$.
\end{prop0}
\proof The points of the projective space $\P d={\bf P}(S^d\R^2)$ correspond to  forms
$f=\sum_{i=0}^d{d\choose i}\,a_i\,x^{d-i}y^i$ , which have coordinates $(a_0,\ldots ,a_d)$ . The rational normal curve $C_d$,
corresponds to polynomials which are $d$-th powers of linear forms. From the expansion $(t_0x+t_1y)^d=\sum_{i=0}^d{d\choose i}\,t_0^{d-i}t_1^i\,x^{d-i}y^i$ we get that the curve $C_d$ can be parametrized by $a_i=t_0^{d-i}t_1^{i}$.
Pick $d-1$ general points on $C_d$ corresponding to $l_i^d=(l_{i,0}x+l_{i,1}y)^d$ for $i=1,\ldots ,d-1$. 
The linear span of $f$ and these points
is a hyperplane, whose equation $\sum {d\choose i}a_ic_i=0$ restricts to $C_d$ to the binary form $\sum {d\choose i}c_it_0^{d-i}t_1^{i}$ of degree $d$ with the $d-1$ real roots
$(t_0,t_1)=(l_{i,0},l_{i,1})$  (because $\sum {d\choose i}c_il_{i,0}^{d-i}l_{i,1}^{i}=0$)
hence also the last root is real, corresponding to a last linear form $l_d^d$. This means that $f$ is a projective linear combination
of the the powers $l_i^d$ for $i=1,\ldots ,d$, or equivalently, $f$ has rank $\le d$.
\qedd

Because of the bijection mentioned earlier, 
$f\in Sym^2(\K^2)$ corresponds to a symmetric $2\times 2$ matrix,
 and its rank is the same as the rank of the associated matrix. Hence $f$ has rank $2$ if and only if 
 $f$ has distinct roots, both for $\K=\R$ or $\K=\C$.
 So $S_{2,r}$ has non empty interior only if $r=2$.
 
For any field of characteristic zero, the dual ring of differential operators $D=\K[\partial_x,\partial_y]$ acts on the polynomial ring
$R=\K[x,y]=\oplus_{d\ge 0} Sym^d(\K^2)=\oplus_{d\ge 0}R_d$ with the usual rules of differentiations. This action preserves the degrees, 
in the sense that we have the pairing
$$R_d\otimes D_k\to R_{d-k}$$
which means that the homogeneous differential operators of degree $k$ takes homogeneous polynomials
of degree $d$ to  homogeneous polynomials
of degree $d-k$.

If $l=ax+by\in R_1$, we denote by $l^{\perp}=b\partial_x-a\partial_y\in D_1$ its orthogonal (apolar) operator,
such that $l^{\perp}(l)=0$.
The main remark, which goes back to the XIXth century, is that the form
 $f=\sum_{j=1}^rc_j(l_j)^d$ is killed by
the operator $\prod_{j=1}^r l_j^{\perp}$.

The space of operators of degree $k$ which kill a given form $f$ of degree $d$
is the kernel of the linear map $A_f\colon D_k\to R_{d-k}$;
see \cite{EhreR93:ejc} for instance.

In the basis $\{\partial_x^{k-i}\partial_y^i\}_{i=0,\ldots k}$
and $\{{d-k\choose j}x^{d-k-j}y^j\}_{j=0,\ldots ,d-k}$ the matrix of $A_f$
is (up to scalar multiples) the well known Hankel matrix (called catalecticant by Sylvester)
of size $(d-k+1)\times (k+1)$

$$\left[\begin{array}{ccccc}a_0&a_1&a_2&\ldots&a_k\\a_1&a_2&a_3&\ldots&\\ \vdots&&&&\vdots\\
a_{d-k} & & & \ldots&a_d\end{array}\right]$$

For general forms $f$, the map $A_f$ has maximal rank. We may assume
$2k\ge d$, so that its kernel has dimension $d-2k$. This holds on any field $\K$.

In particular the general form of degree $d=2n+1$ is killed by a unique operator of degree $n$.
This gives a proof of Sylvester Theorem over $\C$.
The kernel of the Hankel matrix of size $n\times (n+1)$ is an operator $g$
of degree $n$, and decomposing it over $\C$ as  $g=\prod_{j=1}^{n} l_j^{\perp}$
gives the decomposition of $f$ as sum of $n$ $d$-th powers of the $l_j$.
If $f$ is real, the kernel $g$ is real, but the decomposition is possible only if all the roots of $g$  (thought as a polynomial) are real.

Comon and Mourrain explored further this situation \cite{ComoM96:SP}.
If it is not possible to find a kernel of degree $n$ such that all its roots are real,
we consider operators of degree $(n+1)$ and we look for such an operator in the kernel such that all the roots are real.
If this is not possible we go further, and so on.

We use this technique by reproving the classically known criterion
to determine the rank in the case $d=3$.
For $f\in Sym^3(\R^2)$, consider the discriminant
 $$\Delta(f)=4\left|\begin{array}{cc}a_0&a_1\\a_1&a_2\end{array}\right|\cdot
 \left|\begin{array}{cc}a_1&a_2\\a_2&a_3\end{array}\right|-
 \left|\begin{array}{cc}a_0&a_2\\a_1&a_3\end{array}\right|^2$$
 \begin{prop0}\label{d3}
 $S_{3,r}^{\R}$ has non empty interior
 only for $r=2,3$. Precisely, let $f$ be a polynomial of third degree without multiple roots.
Then

(i) $f$ has rank two if and only if $\Delta(f)<0$, or equivalently, if and only if $f$ has one real root

(ii) $f$ has rank three if and only if $\Delta(f)>0$, or equivalently, if and only if $f$ has three real roots.
\end{prop0}
\proof

The differential operators of degree two which annihilate $f$ consist of the kernel of the matrix

$$\left[\begin{array}{ccc}a_0&a_1&a_2\\
a_1&a_2&a_3\\
\end{array}\right]$$

The discriminant of the quadratic generator of the kernel coincides with $-\Delta(f)$; thus 
the operators have two real roots if $\Delta(f)<0$ and this means that the rank-2 complex decomposition is actually
real. Note also that a cubic of real rank two can have only one real root.
Indeed the equation $l_1^3+l_2^3=0$ reduces to the three linear equations
$l_1-e^{\frac{n\pi i}{3}}l_2=0$ for $n=0,1,2$. This proves (i). If $\Delta(f)>0$, the quadratic generator has no real root and (ii) follows from \propref{maxrank}.
\qedd
  
{\bf Remark}~ 
The complete description of rank of tensors of dimension $2$ and order $3$ can be found in 
 \cite{GelfKZ94}(ch. 14, ex. 4.5) in the complex case and in \cite{DesiL08:simax} in the real case.

\section{Forms with all real roots, with an application to quartics}
\begin{prop0}\label{dreal}
Let $f$ be a real binary form of degree $d$ with $d$ real distinct roots. Then the (real) rank of $f$ is equal to $d$.
\end{prop0}
\proof
The proof is by induction on $d$. By the previous section, let $d\ge 3$.
Assume the rank is $\le d-1$. Then we get $f=\sum_{i=1}^{d-1}l_i^d$
and since there are infinitely many such decompositions we may assume that $l_{d-1}$ does not divide $f$. Consider the rational function
$$F=\frac{f}{l_{d-1}}$$

Under a linear (real) change of projective coordinates $\phi(x,y)=(x',y')$ with $y'=l_{d-1}$ we get
$G(x',y')=F(\phi^{-1}(x',y'))=\frac{f(\phi^{-1}(x',y'))}{y'^d}$

Then the polynomial $G(x',1)=\sum_{i=1}^{d-2}n_i(x')^d+1$ has $d$ distinct real roots since $f$ had,
(where $\deg n_i=1$) and its derivative
$\frac{d}{dx'}G(x',1)=\sum_{i=1}^{d-2}dn_i(x')^{d-1}\frac{d}{dx'}(n_i(x'))$ has $d-1$ distinct real roots. 
Now $\frac{d}{dx'}G(x',1)$ has  rank $\le d-2$, indeed $\frac{d}{dx'}n_i(x')$ are constants. This contradicts the inductive assumption. 
Hence the assumption was false and the rank of $f$ must exceed $d-1$. The rank of $f$ must eventually be equal to $d$ from Proposition \ref{maxrank}. 
\qedd

\begin{lemma0}\label{Mobius}
It is easy to check that under a general M\"obius transformation $x\mapsto\frac{ax+b}{cx+d}$ with $ad-bc\neq 0$
both  the number of real roots and the rank of a binary form are invariant.
\end{lemma0}

\begin{lemma0}\label{canonical4}
The following are canonical forms under the action of the M\"obius transformation group
for general forms.

\noindent Case $d=4$:
\begin{itemize}
\item $(x^2+y^2)(x^2+ay^2)\quad$ ($a>0$ no real root, $a<0$ two real roots)
\item  $(x^2-y^2)(x^2+ay^2)\quad$ ($a<0$ four real roots)
\end{itemize}
\medskip
\noindent Case $d=5$:
\begin{itemize}
\item $x(x^2+y^2)(x^2+2axy+by^2)\quad$ ($b-a^2>0$ one real root, $b-a^2<0$ three real roots)
\item  $x(x^2-y^2)(x^2+2axy+by^2)\quad$ ($b-a^2<0$ five real roots)
\end{itemize}

\end{lemma0}

\proof We prove just the first case for $d=4$, the other ones being analogous.
When there are two pairs of conjugate roots, they lie in the complex plane on a circle with real center, then a convenient circle inversion makes the 
four roots on a vertical line.
A translation and a homothety  centered at zero conclude the argument.

When there is one pair of conjugate roots, assume that they are $\pm \sqrt{-1}$.
Then consider the tranformations
$x\mapsto\frac{x+c}{-cx+1}$, which preserve $\pm \sqrt{-1}$ and it is easy to show
that a convenient choice of $c$ makes the sum of the other two roots equal to zero.
\qedd

\begin{prop0}\label{polarprop}
 Let $f$ a real binary form of degree $d$ with distinct roots.

(i) If $f$ has $d$ real roots
then for every $(a,b)\neq (0,0)$ the binary form $af_x+bf_y$ has $d-1$ real roots.

(ii) Conversely, if for every $(a,b)\neq (0,0)$ the binary form $af_x+bf_y$ has $d-1$ real roots
and $3\le d\le 5$, then $f$ has $d$ real roots.
\end{prop0}

\proof
(i) Consider that for any substitution $x=at+c,y=bt+d$ with $ad-bc\neq 0$ we have that
$F(t)=f(at+c,bt+d)$ has $d$ real roots, then
$\frac{d}{dt}f(at+c,bt+d)=af_x+bf_y$ has $d-1$ real roots corresponding to the $d-1$ extremal points of $F$.

(ii) Assume that $f$ has $\le d-1$ (hence $\le d-2$) real roots, and let us show that there exist $(a,b)$ such that $af_x+bf_y$ has $\le d-2$ (hence $\le d-3$) real roots.

For $d=3$, after a M\"obius transformation,
we may assume that $f=x^3+3xy^2$. Then $f_x=3(x^2+y^2)$ has no real roots.

For $d=4$ we may assume by the \lemref{canonical4} that $f=(x^2+y^2)(x^2+ay^2)$. 
For $a>-1$ we consider $f_x=x\cdot(4x^2+2(a+1)y^2)$ which has only one real root.
For $a<-1$ we consider $f_y=y\cdot(4ay^2+2(a+1)x^2)$
which has only one real root.
For $a=-1$ then $f_x-f_y$ has only one real root.

For $d=5$ we may assume by the \lemref{canonical4} that $f=x(x^2+y^2)(x^2+2axy+by^2)$.
The discriminant of $f_x$ is  (up to a positive scalar multiple)
$D(a,b):=
-540a^2-1584a^4+830b^3-180b^4-180b^2-8192a^6+405b^5+405b-7476a^2b^2+1548a^2b+14784a^4b-396a^2b^3+576a^4b^2-
432b^4a^2$. 
  It can be shown that $f_x$ has zero real roots if $D(a,b)>0$ and two real roots if
$D(a,b)<0$. This concludes the proof. \qedd

\begin{coro0}\label{drank}
Let $f$ be a real binary form of degree $d$ with distinct roots. If $f$ has rank $d$ and $3\le d\le 5$ then it has $d$ real roots.
\end{coro0}
\proof
The proof is by induction on $d$. For $d=3$ it follows from the \propref{d3}. Let $4\le d\le 5$. If $f$ has $\le d-2$ real roots then by 
\propref{polarprop}(ii), there exists $(a,b)\neq (0,0)$ such that the binary form $af_x+bf_y$ has $\le d-3$ real roots. 
Then by the inductive assumption
 $af_x+bf_y$ has rank $\le d-2$. So we get $af_x+bf_y=\sum_{i=1}^{d-2}l_i^{d-1}$.
 Choose $c,d $ such that $ad-bc\neq 0$.
 Let $F(t)=f(at+c,bt+d)$.
 We get that $F'(t)=\sum_{i=1}^{d-2}n_i(t)^{d-1}$ for some degree one polynomials $n_i$ and  by integration there is a constant $K$ and degree one 
polynomials $m_i$ such that $\frac{F(t)}{(bt+d)^d}=\sum_{i=1}^{d-2}\frac{m_i(t)^{d}}{(bt+d)^d}+\frac{K}{(bt+d)^d}$.
 
 With the substitution $t=\frac{dx-yc}{-bx+ay}$ we get that the rank of $f$ is $\le d-1$, which is against the assumption.\qedd
\medskip

{\bf Remark} We do not know if the part (ii) of the \propref{polarprop} holds for $d\ge 6$. If this is true, also the 
\corref{drank} can be extended.

\medskip

  Now we can give our promised application to quartics.

\emph{Proof of the Main Theorem in the case $d=4$}
The proof is a direct consequence of \propref{dreal} and \corref{drank}.\qedd

We end this section with a last general result about forms such that their complex rank is smaller than the general one.

\begin{prop0}\label{realcomplex}
Let $f$ be a real binary form of degree $d$ and complex rank $k$.
Assume that $2\le k < \lfloor\frac{d}{2}\rfloor +1$. Then there only two possibilities:

(i) The real rank is equal to the complex rank

(ii)
The real rank of $f$ is $\ge d-k+2$ (when $k=2$ the equality holds).
\end{prop0}
\proof

Assume that (i) does not hold. This means that the contraction

$$D_k\to R_{d-k}$$
has rank $k$ and that the one dimensional kernel is generated by one operator with at least
two complex conjugate roots.
It follows that also the transpose operator
$$D_{d-k}\to R_{k}$$
has rank $k$, and the operators in the kernel are given exactly by the previous operator times every operator
of degree $d-2k$. In particular no operator in the kernel has all real roots.

This argument works also for the next contraction

$$D_{d-k+1}\to R_{k-1}$$
which has again rank $k$. At the next step it is possible to find an operator in the kernel with all real roots. This concludes the proof.
When $k=2$ the equality holds by \propref{maxrank}.
\qedd

\section{Quintics}

It is relatively easy to check if a real quintic form has rank three.
Indeed we have to compute the kernel $[\beta_0,\beta_1,\beta_2,\beta_3]$ of the matrix

$$\left[\begin{array}{cccc}\alpha_0&\alpha_1&\alpha_2&\alpha_3\\
\alpha_1&\alpha_2&\alpha_3&\alpha_4\\
\alpha_2&\alpha_3&\alpha_4&\alpha_5
\end{array}\right]$$

The coefficients $\beta_i$ are given by the $3\times 3$ minors  of the above matrix.
Now in order to check if $\beta=\sum_{i=0}^3x^{3-i}y^i\beta_i$ has three real roots
we consider the invariant of degree $12$ 
$$I_{12}(\alpha_0,\ldots ,\alpha_5):=\Delta(\beta)$$
where $\Delta$ is the discriminant considered before.
A nickname for $I_{12}$ could be the ``apple invariant'', as we will se in the next pictures.
We have that the general $f$ has rank three if and only if $\beta$ has three real roots and by \propref{d3} (ii)
this happens if and only if $I_{12}>0$.

\emph{Proof of the Main Theorem in the case $d=5$}

By \propref{dreal} and \corref{drank} we get that, for $d=5$, a form $f$ with distinct roots has rank $5$
if and only if it has $5$ real roots.
Hence, we have just to prove that both rank three and four have nonempty interior in the other cases.

For quintics with less than 5 real roots, it is enough to consider  the canonical form
$f=x(x^2+y^2)(x^2+2axy+by^2)$ of \lemref{canonical4},
and compute for this form the sign of the invariant $I_{12}$.
Now $f$ has rank three if and only if $I_{12}>0$,
and rank four if and only if $I_{12}<0$.

This can be checked explicitly from the expression (up to a positive scalar multiple)
$I_{12}(a,b)=
2b^3-142b^4-64a^6+3752b^5-83a^2b^2+144a^4b+3960a^2b^3-11840a^4b^2-62900b^4a^2-311952b^5a^4-
593208b^7a^2+375552b^3a^6+287488ba^8+348552b^5a^2+243056a^4b^3+12096a^6b-243410b^6a^2-
673952a^4b^4-392736a^6b^2+99840a^8b^2+700160a^6b^4+700160a^4b^6+92940a^2b^8+7680a^8b^4+
19712a^{10}b^2-336640a^8b^3-6560a^6b^6-287040a^6b^5-115712a^{10}b+480a^4b^8-44400a^4b^7-
4680a^2b^9+77a^2b^{10}-4096a^8-43192b^9+170652b^7-43192b^6+170652b^8-78848a^{10}+3752b^{10}+2b^{12}-
142b^{11}+8192a^{12} $.

The results of the computations are shown in the following figures, and complete the proof.

\begin{figure}
    \centerline{\includegraphics[width=110mm]{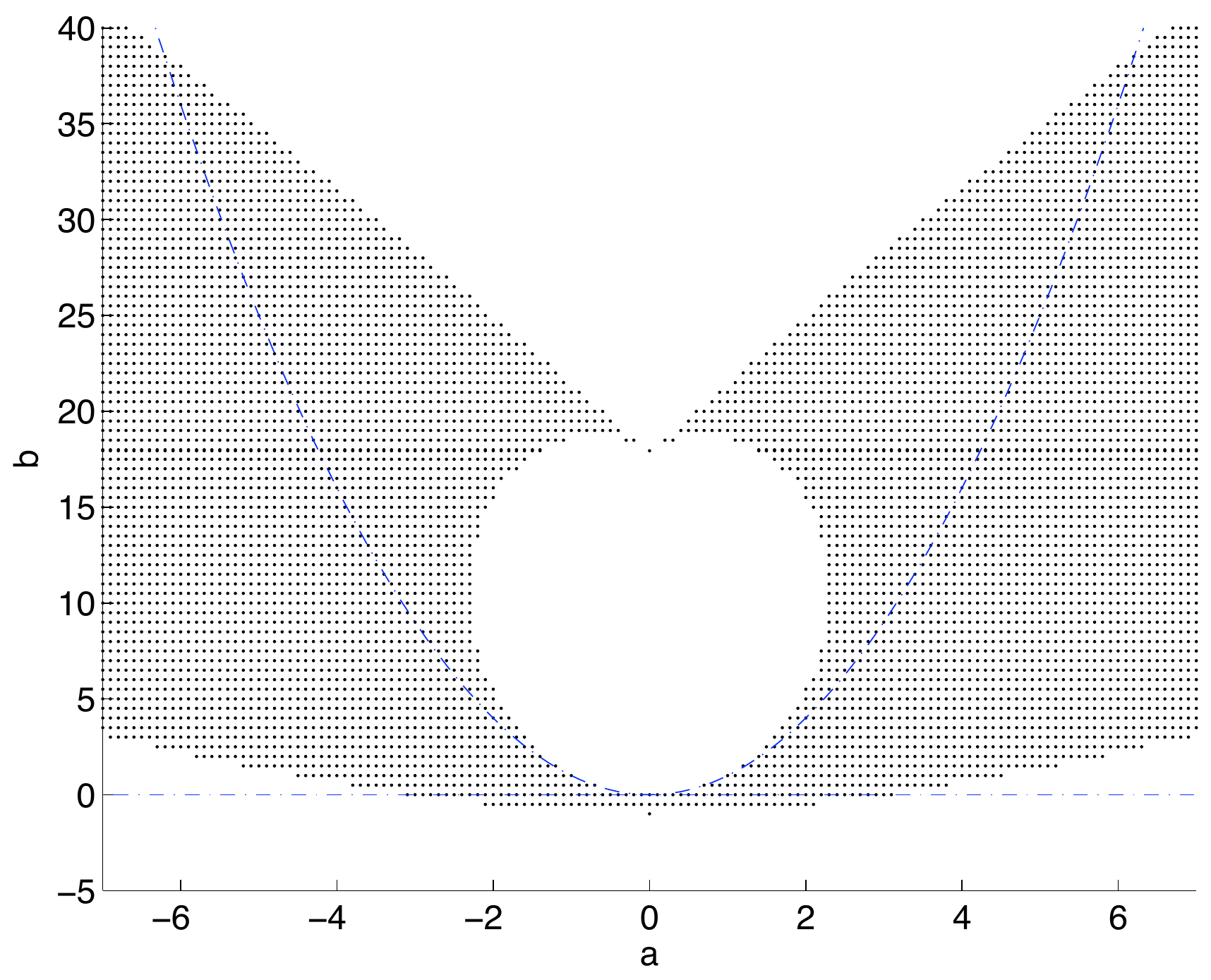}}
\caption{In the shaded area, $I_{12}<0$, and the quintics considered have rank 4. In the area where $I_{12}>0$, they have rank 3.}\label{applePlot-fig}
\end{figure}

For real quintics in the above canonical form the ones with rank four consist of the shaded area in the $(a,b)$-plane in Figure \ref{applePlot-fig}, that is, the area where $I_{12}<0$.
$I_{12}$ contains only even powers of $a$, hence the pictures are symmetric with respect the $b$-axis.

\begin{figure}
    \centerline{\includegraphics[width=110mm]{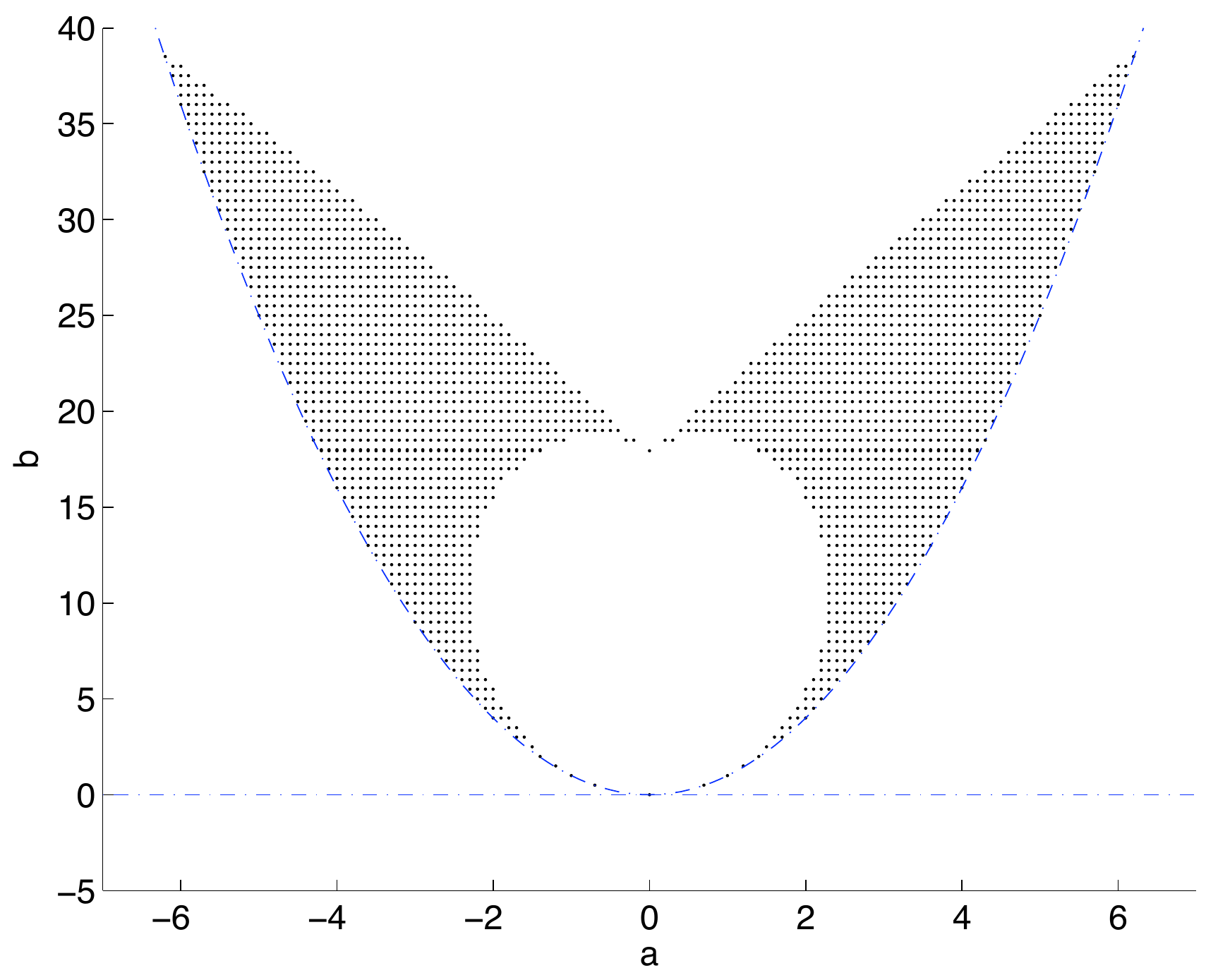}}
    \caption{In the shaded area, the quintics considered have rank four and have only one real root; the parabola $b-a^2=0$ is represented in semi-continuous line.}\label{rankFourOneRoot-fig}
    \end{figure}

It is interesting to compare the previous pictures with the parabola $b-a^2=0$.
It follows that the ``apple'' is all over the parabola (depicted in semi-continuous line in the figures).
The figure \ref{rankFourOneRoot-fig} shows the quintics with one real root (over the parabola) which have rank four; we get a bounded region.

We get that most of quintics with one real root have rank three, but both rank three and four have nonempty interior,
although the rank four part is bounded.
In the same way most of quintics with three real roots (under the parabola) have rank three, but both ranks three and four have nonempty interior,
and even both regions are unbounded.

We conclude this paper with an algorithm which summarizes our results on the quintic.

\subsection*{Algorithm to compute the rank of a generic quintic}

Let
$$f(x)=\alpha_0x^5+5\alpha_1x^4+10\alpha_2x^3+10\alpha_3x^2+5\alpha_4x+\alpha_5$$
\begin{enumerate}
\item
If $f$ has five (distinct) real roots then the real rank of $f$ is $5$.
\item
Otherwise consider the cubic

$$g(x)=\left|\begin{array}{ccc}\alpha_1&\alpha_2&\alpha_3\\\alpha_2&\alpha_3&\alpha_4\\\alpha_3&\alpha_4&\alpha_5\end{array}\right|x^3-
\left|\begin{array}{ccc}\alpha_0&\alpha_2&\alpha_3\\\alpha_1&\alpha_3&\alpha_4\\\alpha_2&\alpha_4&\alpha_5\end{array}\right|x^2+
\left|\begin{array}{ccc}\alpha_0&\alpha_1&\alpha_3\\\alpha_1&\alpha_2&\alpha_4\\\alpha_2&\alpha_3&\alpha_5\end{array}\right|x-
\left|\begin{array}{ccc}\alpha_0&\alpha_1&\alpha_2\\\alpha_1&\alpha_2&\alpha_3\\\alpha_2&\alpha_3&\alpha_4\end{array}\right|$$
\begin{enumerate}
\item
If $g$ has three (distinct) real roots then the real rank of $f$ is $3$.
\item
If $g$ has one real root and a pair of conjugate (distinct) roots then the real rank of $f$ is $4$.
\end{enumerate}
\end{enumerate}

%%%%%%%%%%%%%%%

%\bibliographystyle{unsrt}
%\bibliography{typicalRanks5}

\bigskip\bigskip

\begin{thebibliography}{10}

\bibitem{SidiBG00:ieeesp}
N.~D. SIDIROPOULOS, R.~BRO, and G.~B. GIANNAKIS.
\newblock Parallel factor analysis in sensor array processing.
\newblock {\em IEEE Trans. Sig. Proc.}, 48(8):2377--2388, August 2000.

\bibitem{Kroo08}
P.~KROONENBERG.
\newblock {\em Applied Multiway Data Analysis}.
\newblock Wiley, 2008.

\bibitem{SmilBG04}
A.~SMILDE, R.~BRO, and P.~GELADI.
\newblock {\em Multi-Way Analysis}.
\newblock Wiley, 2004.

\bibitem{BadeK06:tms}
B.~W. BADER and T.~KOLDA.
\newblock Algorithm 862: {M}atlab tensor classes for fast algorithm
  prototyping.
\newblock {\em ACM Trans. Math. Soft.}, 32:635--653, 2006.

\bibitem{MoruHHPA06:ni}
M.~MORUP and L.~K.~HANSEN et~alterae.
\newblock Parallel factor analysis as an exploratory tool for wavelet
  transformed event-related {EEG}.
\newblock {\em NeuroImage}, 29:938--947, 2006.

\bibitem{CarrC70:psy}
J.~D. CARROLL and J.~J. CHANG.
\newblock Analysis of individual differences in multidimensional scaling via
  n-way generalization of {E}ckart-{Y}oung decomposition.
\newblock {\em Psychometrika}, 35(3):283--319, September 1970.

\bibitem{Como02:oxford}
P.~COMON.
\newblock Tensor decompositions, state of the art and applications.
\newblock In J.~G. McWhirter and I.~K. Proudler, editors, {\em Mathematics in
  Signal Processing {V}}, pages 1--24. Clarendon Press, Oxford, UK, 2002.

\bibitem{HackKT05:jnm}
W.~HACKBUSCH, B.~N. KHOROMSKIJ, and E.~E. TYRTYSHNIKOV.
\newblock Hierarchical {K}ronecker tensor-product approximation.
\newblock {\em Jour. Numer. Math.}, 13:119--156, 2005.

\bibitem{IarrK99}
A.~IARROBINO and V.~KANEV.
\newblock {\em Power sums, Gorenstein algebras and determinantal loci}.
\newblock Springer, 1999.

\bibitem{BurgCS97}
P.~B\"URGISSER, M.~CLAUSEN, and M.~A. SHOKROLLAHI.
\newblock {\em Algebraic Complexity Theory}, volume 315.
\newblock Springer, 1997.

\bibitem{Lick85:laa}
T.~LICKTEIG.
\newblock Typical tensorial rank.
\newblock {\em Linear Algebra Appl.}, 69:95--120, 1985.

\bibitem{lanteit}
J.M.LANDSBERG and Z.TEITLER.
\newblock On the ranks and border ranks of symmetric tensors.
\newblock {\em arXiv:0901.0487}, 2009.

\bibitem{BGI}
A.~BERNARDI, A.~GIMIGLIANO, and M.~ID\`A.
\newblock On the stratification of secant varieties of {V}eronese varieties via
  symmetric rank.
\newblock {\em arXiv:0908.1651}.

\bibitem{BL}
A.~BUCZYNSKI and J.M. LANDSBERG.
\newblock Ranks of tensors and a generalization of secant varieties.
\newblock {\em arXiv:0909.4262}.

\bibitem{cilisurvey}
C.~CILIBERTO.
\newblock Geometric aspects of polynomial interpolation in more variables and
  of {W}aring's problem.
\newblock In {\em European {C}ongress of {M}athematics, {V}ol. {I}
  ({B}arcelona, 2000)}, volume 201 of {\em Progr. Math.}, pages 289--316.
  Birkh\"auser, Basel, 2001.

\bibitem{AH}
J.~ALEXANDER and A.~HIRSCHOWITZ.
\newblock Polynomial interpolation in several variables.
\newblock {\em J. Algebraic Geom.}, 4(2):201--222, 1995.

\bibitem{BO}
M.C. BRAMBILLA and G.~OTTAVIANI.
\newblock On the {A}lexander-{H}irschowitz theorem.
\newblock {\em J. Pure Appl. Algebra}, 212(5):1229--1251, 2008.

\bibitem{Tenb00:psy}
J.~M.~F. ten BERGE.
\newblock The typical rank of tall three-way arrays.
\newblock {\em Psychometrika}, 65(5):525--532, September 2000.

\bibitem{AOP}
H.~ABO, G.~OTTAVIANI, and C.~PETERSON.
\newblock Induction for secant varieties of {S}egre varieties.
\newblock {\em Trans. Amer. Math. Soc.}, 361(2):767--792, 2009.

\bibitem{CGG}
M.V. CATALISANO, A.V. GERAMITA, and A.~GIMIGLIANO.
\newblock Secant varieties of {${\bf{P}}\sp 1\times\ldots\times{\bf{P}}\sp 1$}
  (n-times) are not defective for $n \geq 5$.
\newblock {\em arXiv:0809.1701, to appear in J. of Alg. Geom.}

\bibitem{Rez}
B.~REZNICK.
\newblock Laws of inertia in higher degree binary forms.
\newblock {\em arXiv:0906.5559}, June 2009.

\bibitem{EhreR93:ejc}
R.~EHRENBORG and G.~C. ROTA.
\newblock Apolarity and canonical forms for homogeneous polynomials.
\newblock {\em European Journal of Combinatorics}, 14:157--181, 1993.

\bibitem{Rezn92:mams}
B.~REZNICK.
\newblock Sums of even powers of real linear forms.
\newblock {\em Memoirs of the AMS}, 96(463):1--155, March 1992.

\bibitem{ComasS01:arxiv}
G.COMAS and M.~SEIGUER.
\newblock {\em On the rank of a binary form}, 2001.
\newblock arXiv:math/0112311v1.

\bibitem{ComoTDC09:laa}
P.~COMON, J.~M.~F. ten BERGE, L.~De LATHAUWER, and J.~CASTAING.
\newblock Generic and typical ranks of multi-way arrays.
\newblock {\em Linear Algebra Appl.}, 430(11--12):2997--3007, June 2009.

\bibitem{Fr}
S.~FRIEDLAND.
\newblock On the generic rank of $3$-tensors.
\newblock {\em arXiv:0805.3777}, May 2008.

\bibitem{ComoM96:SP}
P.~COMON and B.~MOURRAIN.
\newblock Decomposition of quantics in sums of powers of linear forms.
\newblock {\em Signal Processing, Elsevier}, 53(2):93--107, September 1996.
\newblock special issue on High-Order Statistics.

\bibitem{GelfKZ94}
I.~M. GELFAND, M.~M. KAPRANOV, and A.~V. ZELEVINSKY.
\newblock {\em Discriminants, Resultants, and Multidimensional Discriminants}.
\newblock Birkh\"auser, 1994.

\bibitem{DesiL08:simax}
V.~De SILVA and L-H. LIM.
\newblock Tensor rank and the ill-posedness of the best low-rank approximation
  problem.
\newblock {\em {SIAM} Journal on Matrix Analysis Appl.}, 30(3):1084--1127,
  2008.

\end{thebibliography}

\bigskip

 \textbf{Acknowledgment}~ 
  This work was begun at the workshop on {\it Geometry and representation theory of tensors for computer science, statistics and other areas},
held in July 2008 at the AIM, California, and presented at 
the Workshop on {\it Tensors and Interpolation} held in Nice in June 2009. We thank the AIM
and J. Landsberg for useful talks. 
The work has been also supported in part by contract ANR-06-Blan-0074 ``Decotes''.

\noindent
\textsc{p. comon} - 
Lab. I3S, Cnrs, University of Nice, BP.121, 2000 route des Lucioles, 06903 Sophia-Antipolis cedex (France).  e-mail: \texttt{pcomon@unice.fr}
 
\medskip
\noindent
  \textsc{g. ottaviani} -
  Dipartimento di Matematica ``U. Dini'', Universit\`a di Firenze, viale Morgagni 67/A, 50134 Firenze (Italy). e-mail: \texttt{ottavian@math.unifi.it}

\end{document}